# Two Triangles with the Same Orthocenter and a Vectorial Proof of Stevanovic's Theorem

Prof. Ion Pătrașcu – The National College "Frații Buzești", Craiova, Romania
Prof. Florentin Smarandache – University of New Mexico, U.S.A.

**Abstract.** In this article we'll emphasize on two triangles and provide a vectorial proof of the fact that these triangles have the same orthocenter. This proof will, further allow us to develop a vectorial proof of the Stevanovic's theorem relative to the orthocenter of the Fuhrmann's triangle.

**Lemma 1**

Let $ABC$ an acute angle triangle, $H$ its orthocenter, and $A', B', C'$ the symmetrical points of $H$ in rapport to the sides $BC, CA, AB$.

We denote by $X, Y, Z$ the symmetrical points of $A, B, C$ in rapport to $B'C', C'A', A'B'$
The orthocenter of the triangle $XYZ$ is $H$.

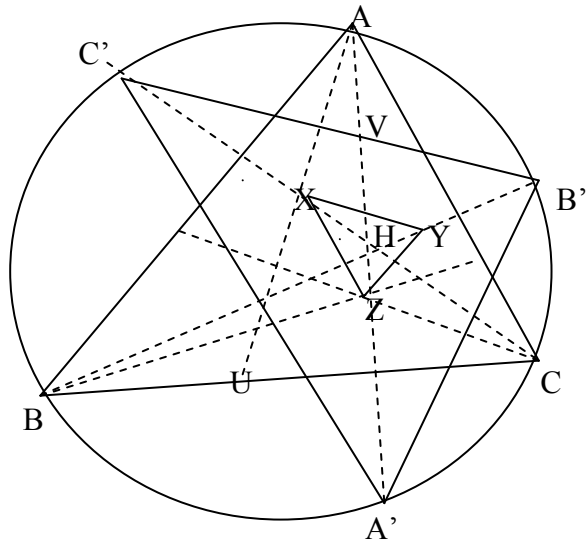

Fig. 1

**Proof**

We will prove that $XH \perp YZ$, by showing that $\overrightarrow{XH} \cdot \overrightarrow{YZ} = 0$.
We have (see Fig.1)

$$\overrightarrow{VH} = \overrightarrow{AH} - \overrightarrow{AX}$$
$$\overrightarrow{BC} = \overrightarrow{BY} + \overrightarrow{YZ} + \overrightarrow{ZC},$$

from here

$$\overrightarrow{YZ} = \overrightarrow{BC} - \overrightarrow{BY} - \overrightarrow{ZC}$$

Because $Y$ is the symmetric of $B$ in rapport to $A'C'$ and $Z$ is the symmetric of $C$ in rapport to $A'B'$, the parallelogram's rule gives us that:

$$\overrightarrow{BY} = \overrightarrow{BC'} + \overrightarrow{BA'}$$
$$\overrightarrow{CZ} = \overrightarrow{CB'} + \overrightarrow{CA'}.$$

Therefore
$$\overrightarrow{YZ} = \overrightarrow{BC} - \left(\overrightarrow{BC'} + \overrightarrow{BA'}\right) + \overrightarrow{B'C} + \overrightarrow{A'C}$$

But
$$\overrightarrow{BC'} = \overrightarrow{BH} + \overrightarrow{HC'}$$
$$\overrightarrow{BA'} = \overrightarrow{BH} + \overrightarrow{HA'}$$
$$\overrightarrow{CB'} = \overrightarrow{CH} + \overrightarrow{HB'}$$
$$\overrightarrow{CA'} = \overrightarrow{CH} + \overrightarrow{HA'}$$

By substituting these relations in the $\overrightarrow{YZ}$, we find:
$$\overrightarrow{YZ} = \overrightarrow{BC} + \overrightarrow{C'B'}$$

We compute
$$\overrightarrow{XH} \cdot \overrightarrow{YZ} = \left(\overrightarrow{AH} - \overrightarrow{AX}\right) \cdot \left(\overrightarrow{BC} + \overrightarrow{C'B'}\right) = \overrightarrow{AX} \cdot \overrightarrow{BC} + \overrightarrow{AH} \cdot \overrightarrow{C'B'} - \overrightarrow{AX} \cdot \overrightarrow{BC} - \overrightarrow{AX} \cdot \overrightarrow{C'B'}$$

Because
$$AH \perp BC$$
we have
$$\overrightarrow{AH} \cdot \overrightarrow{BC} = 0,$$
also
$$AX \perp B'C'$$
and therefore
$$\overrightarrow{AX} \cdot \overrightarrow{B'C'} = 0.$$

We need to prove also that
$$\overrightarrow{XH} \cdot \overrightarrow{YZ} = \overrightarrow{AH} \cdot \overrightarrow{C'B'} - \overrightarrow{AX} \cdot \overrightarrow{BC}$$

We note:
$$\{U\} = AX \cap BC \text{ and } \{V\} = AH \cap B'C'$$
$$\overrightarrow{AX} \cdot \overrightarrow{BC} = AX \cdot BC \cdot cox\sphericalangle(AX, BC) = AX \cdot BC \cdot cox(\sphericalangle AUC)$$
$$\overrightarrow{AH} \cdot \overrightarrow{C'B'} = AH \cdot C'B' \cdot cox\sphericalangle(AH, C'B') = AH \cdot C'A' \cdot cox(\sphericalangle AVC')$$

We observe that
$\sphericalangle AUC \equiv \sphericalangle AVC'$ (angles with the sides respectively perpendicular).

The point $B'$ is the symmetric of $H$ in rapport to $AC$, consequently
$$\sphericalangle HAC \equiv \sphericalangle CAB',$$
also the point $C'$ is the symmetric of the point $H$ in rapport to $AB$, and therefore
$$\sphericalangle HAB \equiv \sphericalangle BAC'.$$

From these last two relations we find that
$$\sphericalangle B'AC' = 2\sphericalangle A.$$

The sinus theorem applied in the triangles $AB'C'$ and $ABC$ gives:
$$B'C' = 2R \cdot \sin 2A$$
$$BC = 2R \sin A$$

We'll show that

$$AX \cdot BC = AH \cdot C'B',$$
and from here
$$AX \cdot 2R \sin A = AH \cdot 2R \cdot \sin 2A$$
which is equivalent to
$$AX = 2AH \cos A$$
We noticed that
$$\sphericalangle B'AC' = 2A,$$
Because
$$AX \perp B'C',$$
it results that
$$\sphericalangle TAB \equiv \sphericalangle A,$$
we noted $\{T\} = AX \cap B'C'$.
On the other side
$$AC' = AH, \quad AT = \frac{1}{2}AY,$$
and
$$AT = AC' \cos A = AH \cos A,$$
therefore
$$\overrightarrow{XH} \cdot \overrightarrow{YZ} = 0.$$
Similarly, we prove that
$$YH \perp XZ,$$
and therefore $H$ is the orthocenter of triangle $XYZ$.

**Lemma 2**

Let $ABC$ a triangle inscribed in a circle, $I$ the intersection of its bisector lines, and $A', B', C'$ the intersections of the circumscribed circle with the bisectors $AI, BI, CI$ respectively.

The orthocenter of the triangle $A'B'C'$ is $I$.

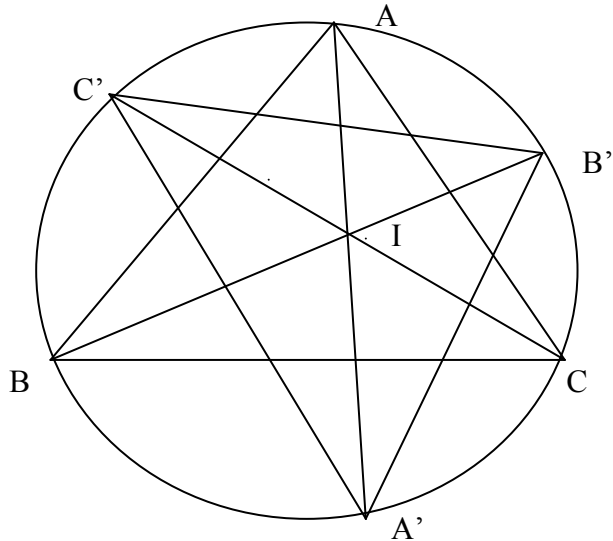

Fig. 2

**Proof**

We'll prove that $A'I \perp B'C'$.

Let
$$\alpha = m\left(\widehat{A'C}\right) = m\left(\widehat{A'B}\right),$$
$$\beta = m\left(\widehat{B'C}\right) = m\left(\widehat{B'A}\right)$$
$$\gamma = m\left(\widehat{C'A}\right) = m\left(\widehat{C'B}\right)$$

Then
$$m\sphericalangle(A'IC') = \frac{1}{2}(\alpha + \beta + \gamma)$$

Because
$$2(\alpha + \beta + \gamma) = 360°$$

it results
$$m\sphericalangle(A'IC') = 90°,$$

therefore
$$A'I \perp B'C'.$$

Similarly, we prove that
$$B'I \perp A'C',$$

and consequently the orthocenter of the triangle $A'B'C'$ is $I$, the center of the circumscribed circle of the triangle $ABC$.

**Definition**

Let $ABC$ a triangle inscribed in a circle with the center in $O$ and $A',B',C'$ the middle of the arcs $\widehat{BC}$, $\widehat{CA}$, $\widehat{AB}$ respectively. The triangle $XYZ$ formed by the symmetric of the points $A',B',C'$ respectively in rapport to $BC, CA, AB$ is called the Fuhrmann triangle of the triangle $ABC$.

**Note**

In 2002 the mathematician Milorad Stevanovic proved the following theorem:

**Theorem** (M. Stevanovic)

In an acute angle triangle the orthocenter of the Fuhrmann's triangle coincides with the center of the circle inscribed in the given triangle.

**Proof**

We note $A'B'C'$ the given triangle and let $A, B, C$ respectively the middle of the arcs $\widehat{B'C'}$, $\widehat{C'A'}$, $\widehat{A'B'}$ (see Fig. 1). The lines $AA', BB', CC'$ being bisectors in the triangle $A'B'C'$ are concurrent in the center of the circle inscribed in this triangle, which will note $H$, and which, in conformity with Lemma 2 is the orthocenter of the triangle $ABC$. Let $XYZ$ the Fuhrmann triangle of the triangle $A'B'C'$, in conformity with Lemma 1, the orthocenter of $XYZ$ coincides with $H$ the orthocenter of $ABC$, therefore with the center of the inscribed circle in the given triangle $A'B'C'$.